\documentclass{article}[12pt]
\usepackage{graphics}

\def\bbbc{{\mathchoice {\setbox0=\hbox{$\displaystyle\rm C$}\hbox{\hbox
to0pt{\kern0.4\wd0\vrule height0.9\ht0\hss}\box0}}
{\setbox0=\hbox{$\textstyle\rm C$}\hbox{\hbox
to0pt{\kern0.4\wd0\vrule height0.9\ht0\hss}\box0}}
{\setbox0=\hbox{$\scriptstyle\rm C$}\hbox{\hbox
to0pt{\kern0.4\wd0\vrule height0.9\ht0\hss}\box0}}
{\setbox0=\hbox{$\scriptscriptstyle\rm C$}\hbox{\hbox
to0pt{\kern0.4\wd0\vrule height0.9\ht0\hss}\box0}}}}
\def\bbbg{{\mathchoice {\setbox0=\hbox{$\displaystyle\rm G$}\hbox{\hbox
to0pt{\kern0.28\wd0\vrule height0.93\ht0\hss}\box0}}
{\setbox0=\hbox{$\textstyle\rm G$}\hbox{\hbox
to0pt{\kern0.28\wd0\vrule height0.93\ht0\hss}\box0}}
{\setbox0=\hbox{$\scriptstyle\rm G$}\hbox{\hbox
to0pt{\kern0.28\wd0\vrule height0.93\ht0\hss}\box0}}
{\setbox0=\hbox{$\scriptscriptstyle\rm G$}\hbox{\hbox
to0pt{\kern0.28\wd0\vrule height0.93\ht0\hss}\box0}}}}
\def\bbbz{{\mathchoice {\hbox{$\mathsf\textstyle Z\kern-0.5em Z$}}
{\hbox{$\mathsf\textstyle Z\kern-0.5em Z$}}
{\hbox{$\mathsf\scriptstyle Z\kern-0.4em Z$}}
{\hbox{$\mathsf\scriptscriptstyle Z\kern-0.3em Z$}}}}
\def\bbbr{{\rm I\!R}} 

\def\bbbn{{\rm I\!N}} 

\def\bbbl{{\rm I\!L}}

\newtheorem{lemma}{Lemma}

\newtheorem{theorem}{Theorem}

\newtheorem{remark}{Remark}

\title{\bf On Emergence of Dominating Cliques in Random Graphs}

\author{Martin Neh\'ez  \\
\small{Department of Information Technologies,} \\
\small{VSM School of Management, City University of Seattle,} \\
\small{Pan\'onska cesta 17,} \\
\small{851 04 Bratislava, Slovak Republic} \\
\small{e-mail: {\tt mnehez@vsm.sk}} \\ \\
Daniel Olej\'ar \\
\small{Department of Computer Science,} \\
\small{FMPI, Comenius University in Bratislava, Mlynsk\'a dolina,} \\
\small{842 48 Bratislava, Slovak Republic} \\ \\
Michal Demetrian \\
\small{Department
of Mathematical and Numerical Analysis,} \\
\small{FMPI, Comenius University in Bratislava, Mlynsk\'a dolina M 105,} \\
\small{842 48 Bratislava, Slovak Republic}
}

\begin{document}
\maketitle

\begin{abstract} Emergence of dominating cliques in Erd\"os-R\'enyi random graph
model ${\bbbg(n,p)}$ is investigated in this paper. It is shown
this phenomenon possesses a phase transition. Namely, we have
argued that, given a constant probability $p$, an $n$-node random
graph $G$ from ${\bbbg(n,p)}$ and for $r= c \log_{1/p} n$ with $1
\leq c \leq 2$, it holds: (1) if $p > 1/2$ then an $r$-node clique
is dominating in $G$ almost surely and, (2) if $p \leq ( 3 -
\sqrt{5})/2$ then an $r$-node clique is not dominating in $G$
almost surely. The remaining range of probability $p$ is discussed
with more attention. A detailed study shows that this problem is
answered by examination of sub-logarithmic growth of $r$ upon $n$.


\end{abstract}

\section{Introduction}

The phase transition phenomenon was originally observed as a
physical effect. In discrete mathematics, it was originally
described by P. Erd\"os and A. R\'enyi in \cite{ER60}. The most
frequently property of graphs which have been studied with
relation to the phase transitions in random graphs is the
connectivity. The recent surveys of known results concerning this
area can be find in Refs. \cite{Bo01} and \cite{JLR00}, Chapter 5.

Our paper deals with another interesting graph problem that is the
emerging of a dominating clique in a random graph. The theory of
dominating cliques in random graphs has several nontrivial
applications in computer science. The most significant ones are:
(1) heuristics in satisfiability search \cite{CGA05} and (2) the
construction of a space-efficient interval routing scheme with a
small additive stretch for almost all and large-scale distributed
systems \cite{NO05isaac}.

\subsection{Preliminaries and terminology}
Given a graph $G=(V, E)$, a set $S \subseteq V$ is said to be a
\emph{dominating set} of $G$ if each node $v \in V$ is either in
$S$ or is adjacent to a node in $S$. The \emph{domination number}
$\gamma(G)$ is the minimum cardinality of a dominating set of $G$.


A \emph{clique} in $G$ is a maximal set of mutually adjacent nodes
of $G$, i.e., it is a maximal complete subgraph of $G$. The
\emph{clique number}, denoted $cl(G)$, is the number of nodes of
clique of $G$. If a subgraph $S$ induced by a dominating set is a
clique in $G$ then $S$ is called a \emph{dominating clique} in
$G$.

The model of random graphs is introduced in the following way. Let
$n$ be a positive integer and let $p \in \bbbr$, $0 \leq p \leq
1$, be a \emph{probability of an edge}. The \emph{(probabilistic)
model of random graphs} $\bbbg(n,p)$ consists of all graphs with
$n$-node set $V = \{ 1, \dots, n \}$ such that each graph has at
most ${n \choose 2}$ edges being inserted independently with
probability $p$. Consequently, if $G$ is a graph with node set $V$
and it has $|E(G)|$ edges, then a probability measure $\Pr$
defined on $\bbbg(n,p)$ is given by:
\[ 
 \Pr [G] = p^{|E(G)|} (1-p)^{{n \choose 2} - |E(G)|}~.
\]
This model is also called \emph{Erd\"os-R\'enyi random graph
model} \cite{{Bo01},{JLR00}}.

Let $A$ be any set of graphs from $\bbbg(n,p)$ with a property
$Q$. We say that \emph{almost all graphs} have the property $Q$
iff:
\[
  \Pr[A] \to 1 \quad \emph{as} \quad n \to \infty~.
\]
The term "almost surely" stands for "with the probability
approaching $1$ as $n \to \infty$".

\subsection{Previous work and our result}

Dominating sets and cliques are basic structures in graphs and
they have been investigated very intensively. To determine whether
the domination number of a graph is at most $r$ is an NP-complete
problem \cite{GJ79}. The maximum-clique problem is one of the
first shown to be NP-hard \cite{Ka72}. A well-known result of B.
Bollob\'as, P. Erd\"os et al. states that the clique number in
random graphs $\bbbg(n,p)$ is bounded by a very tight bounds
\cite{{Bo01},{BE76},{Klb72},{Ma76},{OT97},{Pa}}. Let $b=1/p$ and
let
\begin{equation}
\label{r0}
 r_0 = \log_b n - 2 \log_b \log_b n + \log_b 2 + \log_b \log_b e ~,
\end{equation}
\begin{equation}
\label{r1}
 r_1 = 2 \log_b n - 2 \log_b\log_b n + 2\log_b e + 1 - 2\log_b 2
 ~.
\end{equation}
J. G. Kalbfleisch and D. W. Matula \cite{{Klb72},{Ma76}} proved
that a random graph from $\bbbg(n, p)$ does not contain cliques of
the order greater than $\lceil r_1 \rceil$ and less or equal than
$\lfloor r_0 \rfloor$ almost surely. (See also
\cite{{BE76},{OT97},{Pa}}.) The domination number of a random
graph have been studied by B. Wieland and A. P. Godbole in
\cite{WG01}.

The phase transition of dominating clique problem in random graphs
was studied independently by M. Neh\'ez and D. Olej\'ar in
\cite{{NO05isaac},{NO05KAM}} and J. C. Culberson, Y. Gao, C. Anton
in \cite{CGA05}. It was shown in \cite{CGA05} that the property of
having a dominating clique is monotone, it has a phase transition
and the corresponding threshold probability is $p^* =
(3-\sqrt{5})/2$. The standard first and the second moment methods
(based on the Markov's and the Chebyshev's inequalities,
respectively, see \cite{{AS2000},{JLR00}}) were used to prove this
result. However, the preliminary result of M. Neh\'ez and D.
Olej\'ar \cite{NO05KAM} pointed out that to complete the behavior
of random graphs in all spectra of $p$ needs a more accurate
analysis, namely in the case when $(3-\sqrt{5})/2 < p \leq 1/2$.
The main result of this paper is the refinement of the previous
results from \cite{{CGA05},{NO05isaac},{NO05KAM}}. Let us
formulate this as the following theorem.

\begin{theorem}
\label{Thm1} Let $0 < p < 1$ be fixed and let $\bbbl x$ denote
$\log_{1/(1-p)} x$. Let $r$ be order of a clique such that
$\lfloor r_0 \rfloor \leq r \leq \lceil r_1 \rceil$. Let
$\delta(n): {\bbbn} \to {\bbbn}$ be an arbitrary slowly increasing
function such that $\delta(n) = o (\log n)$ and let $G \in
{\bbbg(n,p)}$ be a random graph. Then:
\begin{enumerate}
\item If $p > 1/2$, then an $r$-node clique
is dominating in $G$ almost surely;
\item If $p \leq ( 3 - \sqrt{5})/2$, then
an $r$-node clique is not dominating in $G$ almost surely;
\item If $(3-\sqrt{5})/2 < p \leq  1/2$, then an $r$-node clique:
\begin{itemize}
\item is dominating in
$G$ almost surely, if $r \geq {\bbbl n} + \delta(n)$,
\item is not dominating in $G$ almost surely, if $r \leq {\bbbl n}  -
\delta(n)$,
\item is dominating with a finite probability $f(p)$ for a suitable function
$f:[0,1] \to [0,1]$, if $r = {\bbbl n} + O(1)$.
\end{itemize}
\end{enumerate}
\end{theorem}
To prove Theorem \ref{Thm1} the first and the second moment method
were used. The leading part of our analysis follows from a
property of a function defined as a ratio of two random variables
which count dominating cliques and all cliques in random graphs,
respectively.
The critical values of $p$: $(3-\sqrt{5})/2$ and $1/2$,
respectively, are obtained from the bounds (\ref{r0}), (\ref{r1})
see \cite{{Klb72},{Ma76}}.

The rest of this paper contains the proof of the Theorem
\ref{Thm1}. Section 2 contains the preliminary results. An
expected number of dominating cliques in ${\bbbg}(n,p)$ is
estimated here. The main result is proved in section 3. Possible
applications are discussed in 
section 4.

\section{Preliminary results}

For $r>1$, let $S$ be an $r$-node subset of an $n$-node graph $G$.
Let $A$ denote the event that "$S$ is a dominating clique of $G
\in {\bbbg}(n,p)$". Let $in_r$ be the associated $0$-$1$
(indicator) random variable on $\bbbg(n,p)$ defined as follows:
$in_r=1$ if $G$ contains a dominating clique $S$ and $in_r=0$,
otherwise. Let $X_r$ be a random variable that denotes the number
of $r$-node dominating cliques. More precisely, $X_r= \sum  in_r$
where the summation ranges over all sets $S$. The following lemma
expresses the expectation of $X_r$.

\begin{lemma} \cite{NO05isaac} The expectation $E(X_r)$ of the random variable $X_r$
is given by:
\begin{equation}
\label{Exr}
 E(X_r) = {n \choose r} p^{{r \choose 2}}(1 - p^r - (1-p)^r)^{n-r} ~.
\end{equation}
\end{lemma}
We use the following properties adopted from \cite{OT97}, pp.
501--502.

\medskip
\noindent{\bf Claim 1.}
\emph{Let $0 < p < 1$ and $k \leq (\eta -1)\frac{\ln n}{\ln p}$,
$\eta<0$ starting with some positive integer $n$. Then:}
\[
(1-p^k)^n =  \exp(-np^k)\left( 1 + O(np^{2k}) \right) =
  1 - np^k + O\left(np^{2k}\right) ~.
\]

\noindent{\bf Claim 2.}
\emph{Let $k = o(\sqrt{n})$, then:}
\[
  n^{\underline{k}} = n(n-1) \cdots (n-k+1)
  = n^k \left( 1 - {k \choose 2} \frac{1}{n}
  + O \left(  \frac{k^4}{n^2} \right) \right)~.
\]
The upper bound on $r$ in $\bbbg(n, p)$ is stated in the following
lemma.
\begin{lemma}
\label{Lemma_ru} Let $b=1/p$ and
\begin{equation}
\label{r_u}
 r_u = 2 \log_b n - 2 \log_b\log_b n + 2\log_b e + 1 - 2\log_b 2
 ~.
\end{equation}
A random graph from $\bbbg(n, p)$ does not contain dominating
cliques of the order greater than $r_u$ with probability
approaching $1$ as $n \to \infty$.
\end{lemma}

\begin{remark}
\label{remark_1} Note that the upper bounds $r_u$ and $r_1$ are
the same.
The argument for estimation of $r_1$ is the same as in Lemma
\ref{Lemma_ru}.
\end{remark}

To obtain conditions for an existence of dominating cliques in
random graphs it is sufficient to estimate the variance
$Var(X_r)$. We can use the fact that the clique number in random
graphs lyes down in a tight interval. We use the bounds (\ref{r0})
and (\ref{r1}) due to \cite{{Klb72},{Ma76}}. The estimation of the
variance $Var(X_r)$ is stated in the following lemma.
\begin{lemma}
\label{var_x} Let p be fixed, $0 < p < 1$ and $\lfloor r_0 \rfloor
\leq r \leq \lceil r_1 \rceil$. Let
\[
 \beta = \min \{~ 2/3, ~-2 \log_b (1-p) ~\}~.
\]
Then:
\begin{equation}
  Var(X_r) = E(X_r) ^2 \cdot O \left( \frac{(\log n)^3}{n^\beta}
  \right)~.
\end{equation}
\end{lemma}

The following claim expresses the number of the dominating cliques
in random graphs.
\begin{lemma}
\label{xr_num} Let $p$, $r$ and $\beta$ be as before, and
\begin{equation}
  X_r = {n \choose r} p^{{r \choose 2}}(1 - p^r - (1-p)^r)^{n-r}
  \times \left\{ 1+ O \left(\frac{(\log n)^{3}}{n^{\beta/2}}
  \right) \right\}  ~.
\end{equation}
The probability that a random graph from $\bbbg(n,p)$ contains
$X_r$ dominating cliques with $r$ nodes is $1-O\left( (\log
n)^{-3} \right)$.
\end{lemma}

\section{Proof of Theorem 1}

For $r>1$, let $Y_r$ be the random variable on $\bbbg(n,p)$ which
denotes the number of $r$-node cliques. According to \cite{OT97},
\begin{equation}
  Y_r = {n \choose r} p^{{r \choose 2}}(1 - p^r)^{n-r}
  \times \left\{ 1+ O \left(\frac{(\log n)^3}{\sqrt{n}}
  \right) \right\}  ~.
\end{equation}
The ratio $X_r/Y_r$ expresses the relative number of dominating
cliques (with $r$ nodes) to all cliques (with $r$ nodes) in
${\bbbg(n,p)}$ and it attains a value in the interval [0, 1]. By
analysis of the asymptotic of $X_r/Y_r$ as $n$ tends $\infty$ we
obtain our main result.

\begin{figure}[h]
\begin{center}
\includegraphics{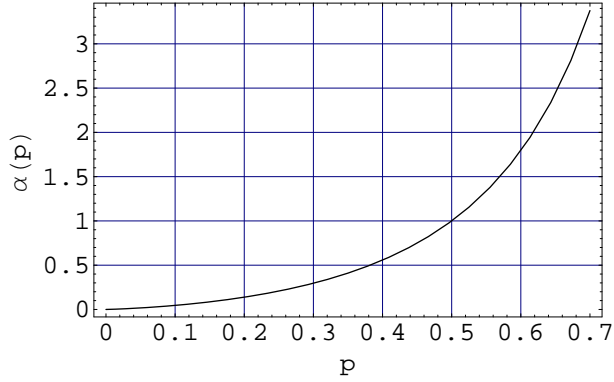}
\caption{The graph of the function $\alpha(p) = - \log_{1/p}
(1-p)$.} \label{Graph_alpha}
\end{center}
\end{figure}

Let us examine the limit value of the ratio $X_r/Y_r$:
$$
 \frac{X_r}{Y_r} = \left(
 \frac{1 - p^r - (1-p)^r}{1 - p^r} \right)^{n-r}
 \times
$$
\begin{equation}
\label{x/y} \times
 \left\{ 1+ O \left(\frac{(\log n)^3}{\sqrt{n}}
  \right) \right\}  \times
 \left\{ 1+ O \left(\frac{(\log n)^3}{n^{\beta/2}}
  \right) \right\}~.
\end{equation}
The most important term of the expression (\ref{x/y}) is the first
one, since the last two terms tend to $1$ as $n \to \infty$. Let
us define $\alpha: [0,1] \to \bbbr$ by:
\[
 \alpha(p) = - \log_{1/p} (1-p)~.
\]
The plot of its graph is in fig. \ref{Graph_alpha} and for the
simplification, we will write also $\alpha$ instead of
$\alpha(p)$. Note that
\begin{equation}
\label{p_alpha_r-x}
 (1-p)^r = p^{r \alpha}~.
\end{equation}

According to Claim 1 and (\ref{p_alpha_r-x}) we have:
\[
 \left(
 \frac{1 - p^r - (1-p)^r}{1 - p^r} \right)^{n-r}
 = \left( 1 - \frac{p^{r \alpha}}{1-p^r}
 \right)^{n-r} =
\]
\[
 = \exp \left(
 \frac{-np^{r \alpha}}{1 - p^r} \right) \cdot
 \left\{ 1+ O \left( n p ^{2r \alpha} \right) \cdot
 \left[ 1 + O\left( \frac{(\log n)^{2+\alpha}}{n}  \right)
 \right]
 \right\} =
\]
\[
 = \exp \left(
 \frac{-np^{r \alpha}}{1 - p^r} \right) \cdot
 \left[
 1 + O \left( n p ^{2r \alpha} \right)
 \right]~.
\]
Let us analyze the asymptotic behavior of the ratio $X_r/Y_r$ as
$n$ tends to $\infty$. According to the assumption $n \to \infty$,
we can write $X_r/Y_r$ in the following two equivalent forms:
\[
 \frac{X_r}{Y_r}=\exp\left( - \frac{np^{r \alpha}}{1 - p^r} \right)
 ~,
\]
or, applying (\ref{p_alpha_r}), as:
\begin{displaymath}
\frac{X_r}{Y_r}=\exp\left( -\frac{n(1-p)^r}{1-p^r}\right) ~.
\end{displaymath}
Using bounds (\ref{r0}) and (\ref{r1}), the admissible number of
nodes of a clique $r$ depends on $n$ as (we consider the leading
term only):
\begin{equation} \label{req0}
r= \rho \log_b n ,
\end{equation}
where $1\leq \rho \leq 2$. This results in:
\begin{displaymath}
\frac{X_r}{Y_r}=\exp\left( -\frac{n^{1-\rho \alpha}}{1-p^r}\right)
,
\end{displaymath}
and one has three different cases:
\begin{enumerate}
\item \ $1-\rho \alpha<0,\ \forall\ \rho \in [1,2]$ \ $\Leftrightarrow$ \
$p>\frac{1}{2}$~,

\item \ $1-\rho \alpha \geq 0,\ \forall\ \rho \in [1,2]$ \ $\Leftrightarrow$ \ $p\leq
\frac{3-\sqrt{5}}{2}$~,

\item \ $1-\rho \alpha$ changes sign as $\rho$ varies in $[1,2]$
\ $\Leftrightarrow $ \ $ \frac{3-\sqrt{5}}{2}<p \leq
\frac{1}{2}$~.
\end{enumerate}
The first case implies
\begin{displaymath}
\lim_{n\to \infty}\frac{X_r}{Y_r}=1 ,
\end{displaymath}
that means the $r-$node cliques is dominating in $G$ almost
surely. The second case implies
\begin{displaymath}
\lim_{n\to \infty}\frac{X_r}{Y_r}=0 ,
\end{displaymath}
and therefore a $r-$node clique is not dominating in $G$ almost
surely. In the third case, there exists a value of $\rho$ (for
each $p$) in the interval $[1,2]$:
\begin{displaymath}
\hat{\rho}=\frac{1}{\alpha(p)} ,
\end{displaymath}
for which we have:
\[
 r = \hat{\rho} \log_b n = \log_{1/(1-p)}n
\]
and
\begin{displaymath}
 \lim_{n\to \infty}\frac{X_r}{Y_r}=
 \exp\left( -n(1-p)^r \right) = e^{-1} .
\end{displaymath}
The ratio $X_r/Y_r$ approaches $1$ ($0$) for $\rho > \hat{\rho}$
($\rho < \hat{\rho}$). Due to corrections of order less than
$\Theta(\log n)$ to the equation (\ref{req0}) taken with
$\rho=\hat{\rho}$ the value of $e^{-1}$ to be changed to another
constant greater or equal than $0$ and less or equal than $1$. The
details are given here. Let $\delta(n): {\bbbn} \to {\bbbn}$ be an
increasing function such that $\delta(n) = o (\log n)$.

If $r=\hat{\rho}\log_{b} n+\delta(n)$, then $X_r/Y_r$ approaches
$1$ as \ $\exp\left(-(1-p)^{\delta(n)}\right)$.

If $r=\hat{\rho}\log_{b} n-\delta(n)$, then $X_r/Y_r$ approaches
$0$ as \ $\exp\left(-(1-p)^{-\delta(n)}\right)$.

And finally, if $r$ differs from $\hat{\rho}\log_{b}n$ by a
constant $\lambda$,  
then the ratio $X_r/Y_r$ asymptotically looks like \
$\exp(-(1-p)^\lambda)$. \\
The proof is complete.  
\ $\diamondsuit$

\begin{figure}[h]
\begin{center}
\includegraphics{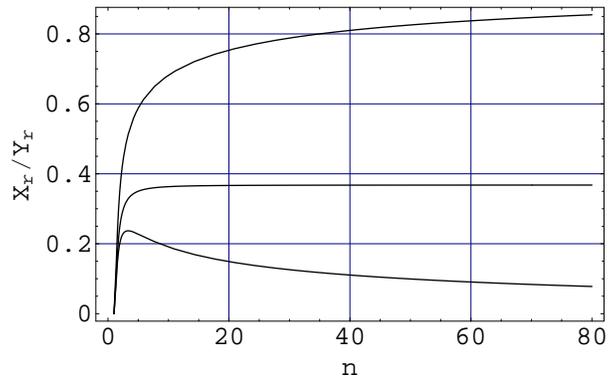}
\caption{The plot of the fraction $X_r/Y_r$ versus $n$ for three
different choices of $\rho$ in the intermediate case when $
\frac{3-\sqrt{5}}{2}<p \leq \frac{1}{2}$. In all three cases $p$
is set to be $0.45$ and $\rho$ varies (from the top to the bottom)
as: $\rho=1.9$, $\rho=1/\alpha(0.45)$, and finally $\rho=1.05$.}
\label{Graph_X/Y}
\end{center}
\end{figure}

\section{Discussion}

We have claimed the conditions for the existence of dominating
cliques in Erd\"os-R\'enyi random graph model. Our result is the
refinement of the previous ones from
\cite{{CGA05},{NO05isaac},{NO05KAM}}.

For possible applications of this result we address the two works
of J. C. Culberson, Y. Gao, C. Anton \cite{CGA05} and M. Neh\'ez
and D. Olej\'ar \cite{NO05isaac}.  The paper \cite{CGA05} deals
with heuristics in satisfiability search. For the second
application, described in \cite{NO05isaac}, we mention the
construction of a space-efficient interval routing scheme with a
small additive stretch in almost all networks modelled by random
graphs $\bbbg(n, p)$ where $p>1/2$. An application of this result
can be found in decentralized content sharing systems based on the
peer-to-peer (shortly P2P) paradigm such as Freenet which uses the
idea of interval routing for retrieving files from local
datastores according to keys \cite{Bon04}.


\medskip
\noindent {\bf Acknowledgement.} This work has been supported by
Gratex Research, Bratislava, by CU grant No. 403/2007 and by the
VEGA grant No. 1/3042/06.

\newpage

\section*{Appendix}

\medskip
\noindent {\bf Proof of Lemma 2.}

The proof follows from the Markov's i\-nequality \cite{JLR00}, p.
8:
\[
 \Pr [~ X \geq t ~] \leq \frac{ E(X) }{ t }~,
 \qquad t > 0 ~.
\label{Markov}
\]
Let us denote $\alpha = \log_{1/p} \left( \frac{1}{1-p} \right) =
- \log_b (1-p)$. Note that:
\begin{equation}
\label{p_alpha_r}
 (1-p)^r = p^{r \alpha}~.
\end{equation}
Let $r=(2 - \varepsilon)\log_b n$, where $0 \leq \varepsilon <1$.
According to Claim 1 we have three cases: $p>1/2,\ p=1/2$ and $p<1/2$. The first two of them can be analyzed together, performing elementary computations
we obtain:
$$
(1-p^r-(1-p)^r)^{n-r}\approx 1-n^{\epsilon-1}\to^{n\to\infty} 1,
~~ \emph{if} \quad p\geq \frac{1}{2} .
$$
In the case $p<1/2$ the same kind of algebra shows that
$$
(1-p^r-(1-p)^r)^{n-r}\approx
\exp\left[-n^{1-\frac{(2-\epsilon)\ln(1-p)}{\ln(p)}}\right], ~~
\emph{if} \quad p<\frac{1}{2} .
$$
We distinguish two different asymptotics in the previous formula. For given $p<1/2$ they are separated by the condition
$$
1-(2-\hat{\epsilon})\frac{\ln(1-p)}{\ln(p)}=0 .
$$
This is solved with respect to $\hat{\epsilon}$ as:
$$
\hat{\epsilon}=2-\frac{\ln(p)}{\ln(1-p)}.
$$
Now we have:
\begin{itemize}
\item for $\epsilon>\hat{\epsilon}$
$$
(1-p^r-(1-p)^r)^{n-r}\to 0 \quad \mbox{as} \quad n\to\infty ,
$$
\item for $\epsilon<\hat{\epsilon}$
$$
(1-p^r-(1-p)^r)^{n-r}\to 1 \quad \mbox{as} \quad n\to\infty ,
$$
\end{itemize}
With respect to upper and lower bound on size of a dominating clique we require $\epsilon$ ranges between $0$ and $1$. This requirement defines
then two critical values of the probability $p$:
\begin{itemize}
\item $\hat{\epsilon}=1$ - in this case
$$
p=\frac{1}{2} ,
$$
\item $\hat{\epsilon}=0$ - in this case
$$
p=\frac{3-\sqrt{5}}{2} .
$$
\end{itemize}
The Stirling's formula (e.g. \cite{Pa}, p. 127) yields
to:
\begin{equation}
\label{compl_sub}
 {n \choose r} p^{{r \choose 2}} \sim \left(
 \frac{n e p^{(r-1)/2} }{r} \right)^r~.
\end{equation}
Consequently,
\[
 {n \choose r_u} p^{{r_u \choose 2}} \to 1 \qquad \emph{and} \qquad
 {n \choose r_u+1} p^{{r_u+1 \choose 2}} \sim \frac{\log_b n}{n} \to 0
\]
The rest follows from the Markov's inequality (\ref{Markov}) for
$t=1$. \ $\diamondsuit$  

\medskip
\noindent {\bf Proof of Lemma 3.}

In order to prove this lemma we will estimate the variance of
$X_r$:
\begin{equation}
\label{var_xe}
  Var(X_r) = E(X^2_r) - E^2(X_r)~.
\end{equation}
The expectation of $X^2_r$ can be expressed in the following way:
$$
E(X^2_r)  = \sum_{j=0}^r {n \choose r} {r \choose j} {n-r \choose
r-j} \cdot {p^2}^{{r \choose 2} - {j \choose 2}} \times
$$
\begin{equation}
\label{e_x2}
 \times ( 1 - p^r - (1-p) ^r ) ^{2n - 4r + 2j} \cdot \Pr
 [ S_r^1,S_r^2 ] ~.
\end{equation}
The equation (\ref{e_x2}) follows from the next analysis. The
nodes of the first dominating clique $S_r^1$ can be chosen in ${n
\choose r}$ ways. The dominating cliques $S_r^1$, $S_r^2$ can (but
need not to) have $j$ common nodes. These nodes can be chosen in
${r \choose j}$ ways. The remaining $(r-1)$ nodes of the second
dominating clique $S_r^2$ have to be chosen from $(n-r)$ nodes of
$V(G) \setminus V(S_r^1 )$. Now we shall choose edges: both
dominating cliques are $r$-node complete graphs and therefore they
contain $2 {r \choose 2}$ edges. But $S_r^1$, $S_r^2$ can have a
nonempty intersection - a complete $j$-node subgraph. Therefore
${j \choose 2}$ edges were counted twice. Both subgraphs $S_r^1$,
$S_r^2$ are dominating cliques and so all $n - 2r + j$ nodes of
the set $V(G) \setminus [ V (S_r^1) \cup V ( S_r^2) ]$ are "good"
with respect to both $S_r^1$, $S_r^2$. The last term,
$\Pr[S_r^1,S_r^2]$ denotes the probability that the nodes of
$V(S_r^1) \setminus V(S_r^2)$ are good with respect to $S_r^2$ and
the nodes of $V(S_r^2) \setminus V(S_r^1)$ are good with respect
to $S_r^1$. It is sufficient to estimate $\Pr [S_r^1, S_r^2 ]$ by
1.

To prove that $Var(X_r)$ is asymptotically less than $E^2(X_r)$,
we extract the expression $E^2(X_r)$ in front of the sum stated by
the equation (\ref{e_x2}). We have:
\begin{equation}
\label{EXr2} E(X_r^2)  \leq E^2(X_r) \cdot \sum_{j=0}^r {n \choose
r}^{-1} {r \choose j} {n-r \choose r-j} \cdot p^{- {j \choose 2}}
\cdot Q(p,r,j)~,
\end{equation}
where $Q(p,r,j) = (1 - p^r - (1 - p)^r)^{-2r + 2j}$.

First we estimate the expression $Q(p,r,j)$. Let us denote $\alpha
= - \log_b (1-p)$, as before. Recall that ~$(1-p)^r = p^{r
\alpha}$. Let us also denote:
\begin{equation}
 \nu = \min \{ 1, - \log_b (1-p) \}~.
\end{equation}

Therefore, from $\lfloor r_0 \rfloor \leq r \leq \lceil r_1
\rceil$ (cf. \cite{OT97}), Claim 1 and (\ref{p_alpha_r}), it
follows:
$$
Q(p,r,j) < \left[ 1- p^{r_0} - p^{\alpha r_0} \right]^{-2r} \leq
$$
$$
\leq \left[ 1 - \frac{(\log_b n)^2}{2n \cdot \log_b e} - \left(
\frac{(\log_b n)^2}{2n \cdot \log_b e } \right)^{\alpha} \right]
^{-4 \log_b n} =
$$
$$
 = \exp  \left\{ 4 \log_b n \cdot
\left[ \frac{(\log_b n)^2}{2n \cdot \log_b e} +
 \left( \frac{(\log_b n)^2}{2n \cdot \log_b e}
 \right)^\alpha \right]
\right\} \times
$$
$$
\times \left( 1 + O \left( \frac {\left( \log n
\right)^{1+2\nu}}{n^{2\nu}}\right) \right) =
$$
$$ = \exp \left(
\frac{2 ( \log_b n )^3}{n \cdot \log_b e} \right) \cdot \exp
\left( \frac {4 ( \log_b n )^{2\alpha+1}}{ (2n \cdot \log_b e
)^{\alpha}}\right)  \cdot \left( 1 + O \left( \frac{ ( \log
n)^{1+2\nu}}{n^{2\nu}} \right) \right),
$$
where $\nu =  \min \{1, \alpha \}$. Since
$$
 \frac{2( \log_b n)^3}{n \cdot \log_b e} \to 0
\qquad  \emph{and} \qquad
 \frac{4( \log_b n)^{2\alpha+1}}
 { ( 2n \cdot\log_b e )^{\alpha}} \to 0
$$
as $n \to \infty$, the value of $Q(p,r,j)$ is $1 + o (1)$ or, more
precisely:
\begin{equation}
\label{Q} Q(p,r,j) = 1 + O \left( \frac{( \log n) ^{2\nu+1}}
{n^{2\nu}}\right)~.
\end{equation}

Now we can concentrate our effort on the estimation of the sum
\begin{equation}
\label{General_sum}
 \sum_{j=0}^r {n \choose  r}^{-1} {r \choose j}
 {n-r \choose r-j} \cdot p^{- {j \choose 2}}~,
\end{equation}
where:
\[
  \lfloor r_0 \rfloor \leq r \leq \lceil r_1 \rceil~.
\]
We use a similar approach as D. Olej\'ar and E. Toman in
\cite{OT97}, pp. 504--506. This sum was also estimated in
Subsection 5.3. of \cite{Pa} (pp. 77--80), but we need more
accurate calculation here. First we introduce the following
notation:
\[
 S(n,r,c,d) = \sum_{j=c}^d {n \choose  r}^{-1} {r \choose j} {n-r \choose r-j}
\cdot b^{{j \choose 2}}~.
\]
Our solution is based on the idea to divide the sum $S(n,r,a,b)$
into three parts by the following way:
\begin{equation}
\label{S_all_ineq}
 S(n,r,0,r) \leq S(n,r,0,1) + S(n,r,2,r_2) + S(n,r,r_2,r)~,
\end{equation}
where:
\[
  r_2 = (1+ \lambda) \log_b n \qquad \emph{for} \qquad 0< \lambda <1~.
\]
All these three parts will be estimated separately. Using Claim 2,
the first part is estimated as follows:
\[
  S(n,r,0,1) = {n-r \choose r} {n \choose  r}^{-1} +
  r \cdot {n-r \choose r-1} {n \choose  r}^{-1} =
\]
\[
  = \left( 1 - \frac{r^2}{n} \right)
  \left[ 1 + O \left( \frac{(\log n)^4}{n^2} \right) \right]
  + \frac{r^2}{n} + O \left( \frac{(\log n)^3}{n^2} \right) =
\]
\begin{equation}
\label{S_0_1}
 = 1 + O \left( \frac{(\log n)^4}{n^2} \right)~.
\end{equation}
To estimate the second part, it is sufficient to analyze the
binomial coe\-fficients. (See also \cite{Pa}, pp. 79--80.)
\[
  {n \choose  r}^{-1} {r \choose j} {n-r \choose r-j} =
  \frac{r!}{n^{\underline{r}}} \cdot
  \frac{r^{\underline{j}}}{j!} \cdot
  \frac{(n-r)^{\underline{r-j}}}{(r-j)!} =
\]
\[
  = \frac{r^{\underline{j}} \cdot (r-j)! }{(r-j)!} \cdot
  \frac{r^{\underline{j}}}{j!} \cdot
  \frac{(n-r)^{\underline{r-j}}}{n^{\underline{j}}
  \cdot (n-j)^{\underline{r-j}}} \ \leq \
  \frac{r^{\underline{j}} \cdot r^{\underline{j}}}{
    j! \cdot n^{\underline{j}} } \ \leq \
  \frac{r^{2j}}{ j! \cdot n^{\underline{j}} }
  \ \sim \ \frac{r^{2j}}{ j! \cdot n^{j} }
\]
We use the Stirling's formula in the following form:
\[
 j! ~\sim~ \left( \frac{j}{e} \right)^j~.
\]
Consequently,
\begin{equation}
  {n \choose  r}^{-1} {r \choose j} {n-r \choose
  r-j} \cdot b^{{j \choose 2}} \ \sim \
  \left(
   \frac{r^2 \cdot b^{j/2} \cdot e}{j \cdot n \cdot \sqrt{b}}
  \right)^j ~.
\end{equation}
The members of the sum $S(n,r,2,r_2)$ attain their asymptotic
maximum for $j=r_2$. More precisely, letting $j=r_2=(1+ \lambda)
\log_b n$ ~we have:
\[
 \frac{r^2 \cdot b^{j/2} \cdot e}{j \cdot n \cdot \sqrt{b}}
 = O \left( \frac{\log n}{n^{1/2-\lambda/2}} \right)~.
\]
Thus,
\[
 S(n,r,2,r_2) \leq \left( \frac{c_1 \cdot \log n}{n^{1/2-\lambda/2}}
 \right)^2 + \left( \frac{c_1 \cdot \log n}{n^{1/2-\lambda/2}} \right)^3
 + \dots + \left( \frac{c_1 \cdot \log n}{n^{1/2-\lambda/2}} \right)^{r_2}~
\]
for a suitable constant $c_1$. It yields:
\begin{equation}
\label{S_2_intermed}
 S(n,r,2,r_2) = O \left( \frac{(\log n)^2}{n^{1-\lambda}}
 \right)~.
\end{equation}

To estimate the sum $S(n,r,r_2,r)$ we extract the term ${n \choose
r}^{-1} \cdot b^{{r \choose 2}}$:
\[
 S(n,r,r_2,r) = {n \choose r}^{-1} \cdot b^{{r \choose 2}}
 \cdot \sum_{j=r_2}^r {r \choose r-j} {n-r \choose r-j}
\cdot p^{{r \choose 2} - {j \choose 2}} ~.
\]
To obtain the upper bound on the right-hand side sum, we
substitute $\lceil r_1 \rceil$ for $r$ in its upper border and
$\lceil r_1 \rceil+1$ for $r$ in all the summands. The reasoning
of such a substitution is the assertion of Lemma \ref{Lemma_ru}
and Remark \ref{remark_1}. We have:
\[
 S(n,r,r_2,r) \leq {n \choose r}^{-1} \cdot b^{{r \choose 2}}
 \cdot \sum_{j=r_2}^{\lceil r_1 \rceil}
 {\lceil r_1 \rceil +1 \choose \lceil r_1 \rceil +1 -j}
 {n- \lceil r_1 \rceil -1 \choose \lceil r_1 \rceil +1 -j}
\cdot p^{{\lceil r_1 \rceil +1 \choose 2} - {j \choose 2}}~.
\]
Let us put $k=\lceil r_1 \rceil+1-j$. Consequently,
\begin{equation}
\label{S_3}
  S(n,r,r_2,r) \leq
\end{equation}
\[
 \leq {n \choose r}^{-1} \cdot b^{{r \choose 2}}
 \cdot \sum_{k=1}^{\lceil r_1 \rceil -r_2+1}
 {\lceil r_1 \rceil +1  \choose k}
 {n- \lceil r_1 \rceil -1  \choose k}
 \cdot p^{k[\lceil r_1 \rceil -(k-1)/2]}~.
\]
Note that
\[
 {\lceil r_1 \rceil +1 \choose k}
 {n- \lceil r_1 \rceil -1 \choose k} \cdot
 p^{k[\lceil r_1 \rceil -(k-1)/2]} \leq
 \left( (\lceil r_1 \rceil +1)\cdot
 np^{\lceil r_1 \rceil -(k-1)/2} \right)^k,
\]
and
\[
 \lceil r_1 \rceil -(k-1)/2 \geq
 \lceil r_1 \rceil/2 + r_2/2 =
\]
\[
 = (3/2 + \lambda/2) \log_b n - \log_b \log_b n +O(1)~.
\]
It yields:
\begin{equation}
\label{S_3rnp}
 ( \lceil r_1 \rceil +1) \cdot
 np^{\lceil r_1 \rceil -(k-1)/2}
 = O \left( \frac{(\log n)^2}{n^{1/2 + \lambda/2}}
 \right)~.
\end{equation}
According to (\ref{S_3}) and (\ref{S_3rnp}),
\[
 S(n,r,r_2,r) \leq {n \choose r}^{-1} \cdot b^{{r \choose 2}}
 \cdot O \left( \frac{(\log n)^2}{n^{1/2 + \lambda/2}}
 \right) ~.
\]
The term ${n \choose r}^{-1} \cdot b^{{r \choose 2}}$ can be
estimated using the Stirling's formula. The estimation is the same
as in the proof of Lemma \ref{Lemma_ru}, see (\ref{compl_sub}).
Thus,
\[
 {n \choose r_1}^{-1} b^{{r_1 \choose 2}} \to 1~,
\]
\[
 {n \choose r}^{-1} b^{{r \choose 2}} \sim
 \frac{(\log_b n)^{c}}{n^c} \to 1~,
\]
if $r=\lceil r_1 \rceil -c$, where $c \geq 1$. Hence,
\begin{equation}
\label{S_3_intermed}
 S(n,r,r_2,r) =
 O \left( \frac{(\log n)^2}{n^{1/2 + \lambda/2}}
 \right)~.
\end{equation}
Let us summarize our results:
\begin{itemize}
\item Eq. (\ref{S_0_1}) shows that $S(n,r,0,1)$ is close to $1$
uniformly with respect to $\lambda$.
\item Eq. (\ref{S_2_intermed}) shows that the "mid" term $S(n,r,2,r_2)$
of the sum-splitting (\ref{S_all_ineq}) is close to zero however,
non-uniformly in $\lambda$. As $\lambda$ approaches $1$ from the
left (i.e. the node number approaches its upper bound)
$S(n,r,2,r_2)$ decreases to zero slowly.
\item Eq. (\ref{S_3_intermed}) shows that $S(n,r,r_2,r)$ is close
to zero uniformly in $\lambda$. (We choose $\lambda=0$ as the
uniform upper bound.)
\end{itemize}

Thus, we have:
$$
E(X^2_r) = E^2(X_r) \cdot \left[ 1 + O \left( \frac{ (\log n)^2
}{n^{2/3}} \right) \right] \cdot \left[ 1 + O \left( \frac {\left(
\log n \right)^{2\nu +1}}{n^{2\nu}}\right) \right]
$$
$$
= E^2(X_r) \cdot \left[ 1 + O \left( \frac {\left( \log n
\right)^{3}}{n^{\beta}}\right) \right]~,
$$
where ~$\nu = \min \{~ 1, - \log_b (1-p) ~\}$~ and \\
$ \beta = \min \{~ 2/3, ~-2 \log_b (1-p) ~\}$~.

Substituting into (\ref{var_xe}) we obtain the estimation of
$Var(X_r)$. \ $\diamondsuit$   

\medskip
\noindent {\bf Proof of Lemma 4.}

It follows from the Chebyshev's inequality \cite{JLR00}: if
$Var(X)$ exists, then:
$$
  \Pr [ | X - E(X) | \geq t ] \geq \frac{ Var(X) }{ t^2 }~,
  \qquad t > 0 ~.
\label{Chebyshev}
$$
Letting $t=E(X_r) \cdot (\log n)^{3} \cdot n^{-\beta /2}$ and
using Lemma \ref{var_x}, we obtain the assertion of Lemma
\ref{xr_num}.  \ $\diamondsuit$  

\end{document}